\documentclass[12pt]{amsart}
\usepackage[square,compress,comma, numbers,sort]{natbib}
\usepackage[colorlinks=true, citecolor=blue, linkcolor=blue]{hyperref}
\usepackage{amsfonts}
\allowdisplaybreaks[4]
\usepackage{graphicx}
\usepackage{epstopdf}

\usepackage{amssymb}
\usepackage{amsmath}
\usepackage{color}
\usepackage{enumerate}
\newcommand{\abs}[1]{\left\lvert #1 \right\rvert}

\def\E#1{\mathbb{E}\left \{#1 \right\}}

\definecolor{c20}{rgb}{0.,0.7,0.}
\definecolor{c30}{rgb}{0.,0.,1.}
\definecolor{c40}{rgb}{1,0.1,0.7}
\definecolor{c50}{rgb}{1,0,0}
\definecolor{c60}{rgb}{1,0.9,0.1}
\definecolor{c70}{rgb}{0.50,1.00,0.00}

\def\N{\mathbb{N}}
\numberwithin{equation}{section}
\newtheorem{theo}{Theorem}[section]
\newtheorem{sat}[theo]{Proposition}
\newtheorem{de}[theo]{Definition}
\newtheorem{lem}{Lemma}[section]

\newtheorem{korr}[theo]{Corollary}
\newtheorem{remark}[theo]{Remark}
\newtheorem{remarks}[theo]{Remarks}
\numberwithin{equation}{section}

\newcommand{\prooftheo}[1]{ \textsc{Proof of Theorem} \ref{#1} }

\newcommand{\prooflem}[1]{\textsc{Proof of Lemma} \ref{#1}}
\newcommand{\proofkorr}[1]{\textsc{Proof of Corollary} \ref{#1}}
\newcommand{\pk}[1]{\mathbb{P} \left\{ #1 \right\} }

\newcommand{\QED}{\hfill $\Box$}
\newcommand{\COM}[1]{}
\def\IF{\infty}
\newcommand{\R}{\mathbb{R}}

\newcommand{\BQN}{\begin{eqnarray}}
\newcommand{\EQN}{\end{eqnarray}}
\newcommand{\BQNY}{\begin{eqnarray*}}
\newcommand{\EQNY}{\end{eqnarray*}}

\def\polhk#1{\setbox0=\hbox{#1}{\ooalign{\hidewidth
\lower1.5ex\hbox{`}\hidewidth\crcr\unhbox0}}}

\newcommand{\kb}[1]{\boldsymbol{#1}}
\newcommand{\vk}[1]{\kb{#1}}

\newcommand{\ve}{\varepsilon}

\usepackage{mathtools}

\newcommand{\norm}[1]{\lVert #1 \rVert}
\newcommand{\nelem}[1]{{Lemma \ref{#1}}}

\newcommand{\nekorr}[1]{{Corollary \ref{#1}}}

\def\vp{\varepsilon}

\def\IF{\infty}

\def\Cov{\mathrm{Cov}}

\date{}

\def\LT{\left}
\def\RT{\right}

\def\Var{\text{Var}}

\newcommand{\limit}[1]{\lim_{#1 \to \infty}}

\newcommand{\BS}{\begin{sat}}
\newcommand{\ES}{\end{sat}}
\newcommand{\BT}{\begin{theo}}
\newcommand{\ET}{\end{theo}}
\newcommand{\BK}{\begin{korr}}
\newcommand{\EK}{\end{korr}}

\newcommand{\BD}{\begin{de}}
\newcommand{\ED}{\end{de}}
\newcommand{\BIT}{\begin{itemize}}
\newcommand{\EIT}{\end{itemize}}
\newcommand{\BDI}{\begin{description}}
\newcommand{\EDI}{\end{description}}
\newcommand{\BRM}{\begin{remarks}}
\newcommand{\ERM}{\end{remarks}}
\newcommand{\BEL}{\begin{lem}}
\newcommand{\EEL}{\end{lem}}

\def\LT{\left}
\def\RT{\right}

\def\Var{\text{Var}}

\def\Cov{\mathrm{Cov}}

\def\LT{\left}
\def\RT{\right}

\def\Var{\text{Var}}

\def\LT{\left}
\def\RT{\right}

\def\Var{\text{Var}}

\def\nj#1{ \mathbb{I}_u \left( #1 \right)}
\def\njk#1{ \mathbb{I}_0 \left( #1 \right)}
\def\njb#1{ \mathbb{I}\left( #1 \right)}
\def\MB{\mathcal{B}}
\def\MF{\mathcal{F}}
\def\Bal{\mathcal{B}_{\alpha}}

  \def\td{\text{\rm d}}

\begin{document}
	
\title{Sojourns of Stationary Gaussian Processes over a Random Interval}
	
	\author{Krzysztof D\c{e}bicki}
	\address{Krzysztof D\c{e}bicki, Mathematical Institute, University of Wroc\l aw, pl. Grunwaldzki 2/4, 50-384 Wroc\l aw, Poland}
	\email{Krzysztof.Debicki@math.uni.wroc.pl}

	\author{Xiaofan Peng}
	\address{Xiaofan Peng, School of Mathematical Sciences, University of Electronic Science and Technology of China, Chengdu 610054, China}
	\email{xfpeng@uestc.edu.cn}

	\bigskip
	
	\date{\today}
	\maketitle

{\bf Abstract}:  We investigate asymptotics of the tail distribution of sojourn time
$$
\int_0^T \mathbb{I}(X(t)> u)dt,
$$
as $u\to\infty$,
where $X$ is a centered stationary Gaussian process and
$T$ is an independent of $X$ nonnegative random variable.
The heaviness of the tail distribution of $T$ impacts the form of the asymptotics,
leading to four scenarios: the case of integrable $T$, the case of regularly varying $T$ with index
$\lambda=1$ and index $\lambda\in(0,1)$  and the case of slowly varying tail distribution of $T$.
The derived findings are illustrated by the analysis of the class of fractional Ornstein-Uhlenbeck processes.



\section{Introduction}
For given stochastic process $X(t),t\ge0$, by
\BQNY\label{defi-sojo}
L_u[a,b] :=\int_{a}^b \nj{X(t)}\td t,
\EQNY
with $\nj{x}:= \njb{x> u}$,
we define
the sojourn time spent above a fixed level $u$ by process $X$ on interval $[a,b]$.
The interest in analysis of distributional properties
of $L_u[a,b]$ stems both from theoretical questions related to the research on
the level sets of stochastic processes and from its importance in applied probability,
as e.g. in finance or insurance theory, where $L_u[0,T]$, $T>0$ may be interpreted as the total time in ruin
up to time $T$ for the risk
process modeled by $X$; see e.g. \cite{R1,R2}.

In the case of $X$ being a Gaussian process, the asymptotics of the tail distribution of $L_u[0,T]$, as $u\to\infty$,
was analyzed extensively in a series of papers by Berman, e.g. \cite{B1,B2}; see also the seminal monograph \cite{Berman92}
and recent refinements \cite{debicki2017approximation,debicki2018sojourn}.

The aim of this paper is to get
the exact asymptotics of tail distribution of $L_u[0,T]$ for a class of centered stationary Gaussian processes
over an independent of $X$ random time $T$.
The motivation to consider extremal behaviour of a stochastic process over a random time interval stems
from its relevance in such problems as ruin of time-changed risk processes \cite{T1,T2}, resetting models \cite{res19}
or hybrid queueing models \cite{ZSD}.
We also refer to related problems on extremes of conditionally Gaussian processes and Gaussian processes with
random variance \cite{Pit1,Pit2}.
Using the fact that
\[
\pk{ L_u[0,T]>0}=\pk{\sup_{t\in[0,T]} X(t)>u},
\]
the findings of this contribution also extend results obtained in \cite{H1,H2,arendarczyk2012exact}.

It appears that the form of the derived exact asymptotics strongly
depends on the heaviness of the tail distribution of $T$, leading to four
scenarios: the case of finite ${\mathbb{E}} T$ (scenario {\bf D1}), the case of $T$ having regularly varying tail distribution
with index $\lambda=1$ (scenario {\bf D2}), $\lambda\in (0,1)$ (scenario {\bf D3}) and the case of
slowly varying tail distribution of $T$ (scenario {\bf D4}); see Section \ref{sect-main-re}.

Brief organisation of the rest of the paper:
In Section \ref{s.notation} we formalize the analyzed model and introduce notation.
In Section \ref{sect-main-re} we derive the tail asymptotic behavior of the sojourn time
for a class of centered stationary Gaussian processes $X$ over random interval $[0,T]$ under introduced
in Section \ref{s.notation} scenarios {\bf D1-D4}, respectively.
Section \ref{s.examples} contains some examples illustrating the main findings of this contribution.
All the proofs are displayed in Section \ref{s.proofs}, whereas few technical results are included in Section \ref{appen}.

\section{Notation and model description}\label{s.notation}

Let $X(t),t\geq0 $ be a centered stationary Gaussian process with a.s. continuous trajectories, unit variance function and covariance function $r$ satisfying
	\begin{itemize}
		\item[{\bf A1}:] $1-r(t)$ is regularly varying at $t=0$ with index $\alpha\in(0,2]$;
        \item[{\bf A2}:]  $r(t)<1$ for all $t>0$;
        \item[{\bf A3}:] $\lim_{t\to\infty}r(t)\log(t)=0$. 
	\end{itemize}
Assumptions {\bf A1-A3} cover wide range of investigated in the literature stationary Gaussian processes, where
{\bf A3} is referred to as {\it Berman's condition} (see, e.g., \cite{Berman92}); see also Section
\ref{s.examples}.

Let function {$v(\cdot)$} be 
such that $\limit{u} v(u)=\IF$  and
 \BQN\label{defi-vu}
\lim_{u\to\IF} u^2 (1-r(1/v(u)))=1.
\EQN
By \cite{Berman92}, {$v(\cdot)$}
exists and is regularly varying at infinity with index $2/\alpha$.

We are interested in the asymptotics of
\[\pk{ L_u^*[0,T]>x}, \]
as $u\to\infty$, where
\BQN \label{defi-sca-sojo}
L_u^*[0,T]:=v(u) L_u[0,T]
\EQN
and $T$ is an independent of $X$ nonnegative random variable with distribution function $F_T(\cdot)$
which belongs to one of the following distribution classes:
\begin{itemize}
      \item[{\bf D1}:] $T$ is integrable;
      \item[{\bf D2}:] $T$ has regularly varying tail distribution with index $\lambda=1$;
      \item[{\bf D3}:] $T$ has regularly varying tail distribution with index $\lambda\in(0,1)$;
       \item[{\bf D4}:] $T$ has slowly varying tail distribution.
     \end{itemize}


Define  for any $x\geq0$
\BQN\label{def-B-H-x}
\MB_{\alpha} (x) =\lim_{S\to\IF} S^{-1}\MB_{\alpha} (S,x),
 \EQN
 with
\BQN\label{BalSx}
 \MB_{\alpha}(S,x)= \int_{\R} \pk{ \int_{0}^{S} \njk{W_\alpha(s) +z } \td s >x } e^{-z} \td z,\quad
 W_\alpha(t)=\sqrt{2}B_{\alpha}(t)-\abs{t}^{\alpha},
 \EQN
where $B_{\alpha}$ is a standard fractional Brownian motion (fBm) with Hurst index $\alpha/2 \in (0,1]$.
By Theorem 2.1 in \cite{debicki2017approximation}, we know that $\MB_\alpha(x)$ is positive and finite for any $x\geq0$. Let $ \mathcal{E}$  be a unit exponential random variable independent of $W_\alpha$ and set
 \BQNY\label{defi-g-ale-x}
\mathcal{G}_{\alpha}(x) = \pk{ \int_{ \R } \njk{W_{\alpha}(s)  +  \mathcal{E}} \td s \leq x }.
\EQNY
As shown in \cite{debicki2018sojourn}, $\mathcal{G}_{\alpha}$ is continuous on $\R^+$, and thus by Remark 2.2 ii) in \cite{debicki2017approximation}
$$\MB_{\alpha}(x)=\int_x^\infty\frac{1}{y}\td \mathcal{G}_{\alpha}(y) $$
holds for all $x\in\R^+$.
We note that $\Bal(0)$ is equal to the classical Pickands constant; see e.g. \cite{Pit96} or Section 10 in \cite{Berman92}.
Let
\BQN\label{def-mu}
m(u)=\big(\Bal(0) v(u) \Psi(u)\big)^{-1},
\EQN
where $\Psi(u)$ is the survival function of an $N(0,1)$ random variable.
Then, by Theorem 10.5.1 in \cite{Berman92}, 
\BQN\label{asym-sup}
\pk{\sup_{t\in[0,1]}X(t)>u}\sim m^{-1}(u),\quad u\to\infty.
\EQN
 In our notation $\sim$ stands for asymptotic equivalence of two
functions as the argument tends to 0 or to $\IF$ respectively.

\section{Main results} \label{sect-main-re}
In this section we find the exact asymptotics of
\begin{eqnarray}
\pk{ L_u^*[0,T]>x}\label{main.1}
\end{eqnarray}
as $u\to\infty$, under scenarios {\bf D1-D4}, respectively.
All the proofs are postponed to Section \ref{s.proofs}.
\subsection{Scenario {\bf D1}}
We begin with the case when $T$ is integrable.
It appears that under this scenario
the main contribution to the asymptotics of (\ref{main.1}) comes from Gaussian process $X$,
whereas $T$ contributes only by its average behavior.
\BT\label{theo-inte}
Let $X(t),t\geq0$ be a centered stationary Gaussian process with unit variance and covariance function satisfying {\bf A1-A2}.
Suppose that $T$ is an independent of $X$ nonnegative random variable that satisfies {\bf D1}. Then for any $x\geq0$
\BQN\label{asym-inte}
\pk{L_u^\ast[0,T]>x} \sim\MB_{\alpha}(x)\E{T}v(u)\Psi(u),\quad u\to\infty.
\EQN
\ET

We can rewrite the result in Theorem \ref{theo-inte} as
\BQNY
\pk{L_u^\ast[0,T]>x} \sim \E{T}\pk{L_u^\ast[0,1]>x},\quad u\to\infty.
\EQNY

\subsection{Scenario D2}
Under this scenario the asymptotics of (\ref{main.1}) is similar to
the one obtained for case {\bf D1} with the exception that
$T$ contributes to (\ref{main.1}) by its integrated tail distribution rather than by
its mean.

\BT\label{theo-lam1}
Let $X(t),t\geq0$ be a centered stationary Gaussian process with unit variance and covariance function satisfying
{\bf A1-A3}. Suppose that $T$ is an independent of $X$ nonnegative random variable that satisfies {\bf D2}.
Then for any $x\geq0$
\BQN\label{asym-lam1}
\pk{L_u^\ast[0,T]>x} \sim 
\MB_{\alpha}(x)l(m(u))v(u)\Psi(u),\quad u\to\infty,
\EQN
 where $l(u)=\int_0^u\pk{T>t}\td t$.
\ET

\begin{remark} \label{remark-rv1}
We note that if $T$ satisfies {\bf D2} and is integrable, then \eqref{asym-lam1} coincides with \eqref{asym-inte}.
\end{remark}

\subsection{Scenario {\bf D3}}
This scenario leads to the asymptotics of (\ref{main.1}) which depends only
of the heaviness of the tail distribution of $T$.

The following continuous distribution function
\BQN \label{def-hf-alp}
\MF_\alpha(x):=\Bal^{-1}(0)\int_0^x\frac{1}{y}\td \mathcal{G}_{\alpha}(x),\quad x\geq0
\EQN
plays an important role in further analysis.
$\overline{\MF_\alpha^{*k}}(x) $ denotes the tail distribution of the $k$-th convolution of
$\MF_\alpha$ at $x\ge0$.

\BT\label{theo-rvso}
Let $X(t),t\geq0$ be a centered stationary Gaussian process with unit variance and covariance function satisfying {\bf A1-A3}.
Suppose that $T$ is an independent of $X$ nonnegative random variable that satisfies {\bf D3}.
 Then for any $x\geq0$
\BQN\label{asym-rvso}
\pk{L_u^\ast[0,T]>x} \sim \lambda\sum_{k=1}^\infty \frac{\Gamma(k-\lambda)}{k!}\overline{\MF_\alpha^{*k}}(x) \pk{T>m(u)},\quad u\to\infty.
\EQN
\ET

\begin{remark} \label{remark-rvso}
Taking $x=0$ in \eqref{asym-rvso} and using
$$\lambda\sum_{k=1}^\infty \frac{\Gamma(k-\lambda)}{k!}=\lambda\sum_{k=1}^\infty \frac{1}{k!} \int_0^\infty l^{k-\lambda-1}e^{-l} \td l =\lambda\int_0^\infty (1-e^{-l})l^{-\lambda-1}\td l=\Gamma(1-\lambda),$$
 we recover Theorem 3.2 in \cite{arendarczyk2012exact}.
\end{remark}

\subsection{Scenario {\bf D4}}
Suppose now that $T$ has slowly varying tail distribution at $\infty$. As shown in the following theorem,
similarly to case {\bf D3},
the asymptotics of (\ref{main.1}) depends only on the asymptotic behavior of the tail distribution
of $T$ but in contrast to scenario {\bf D3} doesn't depend on $x$.
\BT\label{theo-sv}
Let $X(t),t\geq0$ be a centered stationary Gaussian process with unit variance and covariance function satisfying {\bf A1-A3}.
Suppose that $T$ is an independent of $X$ nonnegative random variable that satisfies {\bf D4}. Then for any $x\geq0$
\BQNY\label{asym-sv}
\pk{L_u^\ast[0,T]>x} \sim \pk{T>m(u)},\quad u\to\infty.
\EQNY
\ET
\section{Examples}\label{s.examples}
In this section we illustrate the results derived in Section \ref{sect-main-re} by two classes of
stationary Gaussian processes: fractional Ornstein-Uhlenbeck processes and increments of fractional Brownian motions.

\subsection{Fractional Ornstein-Uhlenbeck processes}
Suppose that $X$ is a centered stationary Gaussian process with covariance
$r(t)=e^{-t^\alpha}, t\geq0$,
for $\alpha\in(0,2]$.
We call $X$ a {\it fractional Urnstein-Uhlenbeck} process with index $\alpha$. If $\alpha=1$, then $X$ is the classical Ornstein-Uhlenbeck process.

It is straightforward to check that {\bf A1-A3} are satisfied.
Thus, the following proposition holds due to Theorems \ref{theo-inte}-\ref{theo-sv}.
\begin{sat}
Suppose that $X$ is a fractional Ornstein-Uhlenbeck process with index $\alpha\in(0,2]$,
and $T$ is an independent of $X$ nonnegative random variable.
Then for any $x\geq0$, as $u\to\IF$,
\\
(i) If $T\in$ {\bf D1}, then $\pk{L_u^\ast[0,T]>x} \sim \MB_{\alpha}(x)\E{T}(2\pi)^{-1/2}u^{2/\alpha-1}e^{-u^2/2}$. \\
(ii) If $T\in$ {\bf D2}, then $\pk{L_u^\ast[0,T]>x} \sim \MB_{\alpha}(x) (2\pi)^{-1/2}u^{2/\alpha-1}e^{-u^2/2} \int_0^{\sqrt{2\pi}\MB_{\alpha}^{-1}(0) u^{1-2/\alpha}e^{u^2/2}}\pk{T>t} \td t.$\\
(iii) If $T\in$ {\bf D3}, then $\pk{L_u^\ast[0,T]>x} \sim \lambda\sum_{k=1}^\infty \frac{\Gamma(k-\lambda)}{k!}
\overline{\MF_\alpha^{*k}}(x) \pk{T>\sqrt{2\pi}\MB_{\alpha}^{-1}(0) u^{1-2/\alpha}e^{u^2/2}}$.\\
(iv) If $T\in$ {\bf D4}, then $\pk{L_u^\ast[0,T]>x} \sim \pk{T>\sqrt{2\pi}\MB_{\alpha}^{-1}(0) u^{1-2/\alpha}e^{u^2/2}}$.
\end{sat}

\subsection{Increments of fractional Brownian motion}
 For a standard fBm $B_\alpha(t), t\geq0$ with Hurst index $\alpha/2\in(0,1)$
and $a>0$, define
 \BQNY
 X_{\alpha,a}(t):=\frac{B_\alpha(t+a)-B_\alpha(t)}{a^{\alpha/2}},\quad t\geq 0.
 \EQNY
One can check that $X_{\alpha,a}$ is a centered stationary Gaussian process with unit variance and covariance function
$$r(t)=\frac{(a+t)^{\alpha}+\abs{a-t}^{\alpha}-2t^{\alpha}}{2a^{\alpha}},\quad t\geq0.$$
and 
$ 1-r(t)
 \sim a^{-\alpha}t^{\alpha},$ $t\to0,
$
which verifies assumption {\bf A1}. Similarly  for $t>a$
 \BQNY
\abs{r(t)}
 \leq \frac{\alpha\abs{1-\alpha}(t-a)^{\alpha-2}}{2a^{\alpha-2}},
  \EQNY
which confirms assumption {\bf A3}.
Thus the following proposition holds.

\begin{sat}
Suppose that $X_{\alpha,a}(t),t\ge 0$ with $\alpha\in(0,2)$, $a>0$ is independent of
a nonnegative random variable $T$.
Then for any $x\geq0$, as $u\to\IF$,
\\
(i) If $T\in$ {\bf D1},
then $\pk{L_u^\ast[0,T]>x} \sim \MB_{\alpha}(x)\E{T}(2\pi)^{-1/2}a^{-1}u^{2/\alpha-1}e^{-u^2/2}$. \\
(ii) If $T\in$ {\bf D2}, then $$\pk{L_u^\ast[0,T]>x} \sim \MB_{\alpha}(x) (2\pi)^{-1/2}a^{-1}u^{2/\alpha-1}e^{-u^2/2} \int_0^{\sqrt{2\pi}a\MB_{\alpha}^{-1}(0) u^{1-2/\alpha}e^{u^2/2}}\pk{T>t} \td t.$$
(iii) If $T\in$ {\bf D3}, then $\pk{L_u^\ast[0,T]>x} \sim \lambda\sum_{k=1}^\infty \frac{\Gamma(k-\lambda)}{k!}\overline{\MF_\alpha^{*k}}(x) \pk{T>\sqrt{2\pi}a\MB_{\alpha}^{-1}(0) u^{1-2/\alpha}e^{u^2/2}};$\\
(iv) If $T\in$ {\bf D4}, then $\pk{L_u^\ast[0,T]>x} \sim \pk{T>\sqrt{2\pi}a\MB_{\alpha}^{-1}(0) u^{1-2/\alpha}e^{u^2/2}}$.
\end{sat}


\section{Proofs}\label{s.proofs}
In this section we give detailed proofs of all the theorems presented in Section \ref{sect-main-re}.
We first give a simple extension of Theorem 7.4.1 of \cite{Berman92}.

\BEL\label{lem-lap-soj}
Let $X(t),t\geq0$ be a centered stationary Gaussian process with unit variance and covariance function satisfying {\bf A1} and {\bf A3}. If $L_u^\ast$, $m(u)$ and $\MF_\alpha$ are defined in \eqref{defi-sca-sojo}, \eqref{def-mu} and \eqref{def-hf-alp}, respectively,
then for any $s\geq0$ and 
$0<l_0<l_1<\infty$ we have
\BQN\label{lap-soj}
\lim_{u\to\infty}\sup_{\tau\in[l_0,l_1]}
\abs{\E{e^{-sL_u^\ast[0,\tau m(u)]}}-e^{-\tau\int_0^\infty(1-e^{-sx})\td \MF_\alpha(x)} }=0.
\EQN
\EEL

\prooflem{lem-lap-soj}
For any $\tau>0$, the point convergence  follows from Berman's proof of Theorem 7.4.1 in \cite{Berman92}.
The uniformity of the convergence on $[l_0,l_1]$ follows by monotonicity
of
$\E{e^{-sL_u^\ast[0,\tau m(u)]}}$
and by continuity of
$e^{-\tau\int_0^\infty(1-e^{-sx})\td \MF_\alpha(x)}$
as function of $\tau$.
\QED

Define a compound Poisson process
\BQN\label{def-cp}
Y(t)=\sum_{i=1}^{N(t)}\xi_i,
\EQN
 where $\{N(t):t\geq 0\}$ is a Poisson process with unit intensity, and $\{\xi_{i}:i\geq 1\}$ are independent and identically distributed random variables, with distribution function $\MF_\alpha$, which are also independent of $N$.
The following corollary of Lemma \ref{lem-lap-soj} will play an important role in the proof of Theorem \ref{theo-rvso}.

 \BK\label{cor-unc}
If $X$ is the Gaussian process given as in \nelem{lem-lap-soj} and $Y$ is defined in \eqref{def-cp}, then for any $x\geq0$ and $0<l_0<l_1<\infty$ we have
 \BQN\label{uniw-conv}
\lim_{u\to\infty}\sup_{l\in[l_0,l_1]}\abs{\pk{L_u^\ast[0,lm(u)]>x}-\pk{Y(l)>x}}=0.
\EQN
 \EK

\proofkorr{cor-unc} For arbitrary $l>0$, by
\eqref{lap-soj}, 
\BQNY
\pk{L_u^\ast[0,lm(u)]>x} \rightarrow \pk{Y(l)>x},\quad u\to\IF
\EQNY
holds for any $x>0$. Further,
\BQNY
\abs{\pk{Y(l_1)>x}-\pk{Y(l_0)>x}}\leq \pk{Y(\abs{l_1-l_0})>0} =1-e^{-\abs{l_1-l_0}}\leq \abs{l_1-l_0},
\EQNY
which implies that for any $x>0$, $\pk{Y(l)>x}$ is continuous in $l$.
Finally, the uniform convergence follows by the same argument as stated in \nelem{lem-lap-soj}.
For $x=0$ in \eqref{uniw-conv}, we refer to Lemma 4.3 in \cite{arendarczyk2012exact}.
This completers the proof. \QED

\subsection{\prooftheo{theo-inte}}
 By \eqref{asym-sup}, for arbitrary $\varepsilon>0$, there exists large enough $u$ such that
$$\pk{\sup_{t\in[0,1]}X(t)>u}< (1+\vp)m^{-1}(u),$$
which together with the stationarity of process $X$ implies that for any $x\geq0$ and $t>0$
\BQNY
\frac{\pk{L_u^\ast[0,t]>x}}{v(u)\Psi(u)} &\leq& \frac{\pk{\sup_{s\in[0,t]}X(s)>u}}{v(u)\Psi(u)}\\
&\leq& (t+1)\frac{\pk{\sup_{s\in[0,1]}X(s)>u}}{v(u)\Psi(u)}\\
&\leq& (t+1)(1+\vp)\MB_{\alpha}(0).
\EQNY
Consequently, for nonnegative random variable $T$ with distribution function $F_T$ satisfying {\bf D1}, by dominated convergence theorem and Remark 2.2 i) in \cite{debicki2017approximation} we have
\BQNY
\lim_{u\to\infty}\frac{\pk{L_u^\ast[0,T]>x}}{v(u)\Psi(u)} &=&
 \lim_{u\to\infty}\int_0^\infty\frac{\pk{L_u^\ast[0,t]>x}}{v(u)\Psi(u)} \td F_T(t)\\
 &=& \MB_{\alpha}(x) \int_0^\infty t \td F_T(t)\\
  &=& \MB_{\alpha}(x) \E{T}.
\EQNY
This completes the proof. \QED

\subsection{\prooftheo{theo-lam1}}
 Let $A(u)$ satisfy
  \BQNY
\lim_{u\to\infty}A(u)v(u)=\IF\quad  \textrm{and}\quad \lim_{u\to\infty}A(u)=0.
\EQNY
 By \nekorr{corsm}, for any $x\geq0$ and arbitrary $\varepsilon\in(0,1)$, there exist $\delta>0$ and $u_0$ such that

\BQNY
\inf_{t\in[A(u),\delta m(u)]}\frac{\pk{L_u^*[0,t]>x}}{t\MB_\alpha(x)v(u)\Psi(u)} \geq 1-2\vp,\quad u>u_0
\EQNY
and
\BQNY
\sup_{t\in[A(u),\delta m(u)]}\frac{\pk{L_u^*[0,t]>x}}{t\MB_\alpha(x)v(u)\Psi(u)} \leq 1+2\vp,\quad u>u_0.
\EQNY
Therefore,
\BQNY
\liminf_{u\to\IF}\frac{\pk{L_u^\ast[0,T]>x}}{\MB_{\alpha}(x)l(m(u))v(u)\Psi(u)}
&\geq& \liminf_{u\to\IF}\frac{\int_{A(u)}^{\delta m(u)}\pk{L_u^\ast[0,t]>x} \td F_T(t)}{\MB_{\alpha}(x)l(m(u))v(u)\Psi(u)}\\
&\geq& (1-2\vp)\liminf_{u\to\IF} \frac{\int_{A(u)}^{\delta m(u)}t \td F_T(t)}{l(m(u))}\\
&=& (1-2\vp)\liminf_{u\to\IF} \frac{\int_{0}^{\delta m(u)}t \td F_T(t)}{l(m(u))}\\
&=& (1-2\vp)\liminf_{u\to\IF} \frac{\int_0^{\delta m(u)}\pk{T>t} \td t -\delta m(u)\pk{T>\delta m(u)}}{l(m(u))}\\
&=& (1-2\vp),
\EQNY
where the last inequality follows by Proposition 1.5.9a in \cite{Bingham1989}
such that $l(u)$ is slowly varying at $\IF$ and $\lim_{u\to\IF}u\pk{T>u}/l(u)=0$.

Similarly,
\BQNY
&&\quad\limsup_{u\to\IF}\frac{\pk{L_u^\ast[0,T]>x}}{\MB_{\alpha}(x)l(m(u))v(u)\Psi(u)}\\
&&\leq \limsup_{u\to\IF}\frac{\pk{L_u^\ast[0,A(u)]>x} \pk{T\leq A(u)}+\int_{A(u)}^{\delta m(u)}\pk{L_u^\ast[0,t]>x} \td F_T(t) +\pk{T>\delta m(u)}}{\MB_{\alpha}(x)l(m(u))v(u)\Psi(u)}\\
&&\leq \limsup_{u\to\IF}\frac{A(u)\pk{T\leq A(u)}}{l(m(u))} + (1+2\vp)\limsup_{u\to\IF} \frac{\int_{A(u)}^{\delta m(u)}t \td F_T(t)}{l(m(u))}\\
&&= (1+2\vp),
\EQNY
where the last inequality follows by \eqref{asysojsmu} and the same reasons as above.
Since $\vp$ is arbitrary, letting $\vp\to0$, we complete the proof.\QED

\subsection{\prooftheo{theo-rvso}}
 First, note that by Raabe's Test, the series in (\ref{asym-rvso})
 converges for $\lambda\in(0,1)$. Then by integration by parts, for any $x\geq0$ we have
 \BQN\nonumber
 \int_0^\infty l^{-\lambda} \td \pk{Y(l)>x} &=&  \int_0^\infty \pk{Y(l)>x} \lambda l^{-\lambda-1} \td l -\lim_{l\to0}l^{-\lambda}\pk{Y(l)>x}\\ \nonumber
 &=& \lambda\sum_{k=1}^\infty \overline{\MF_\alpha^{*k}}(x)\frac{1}{k!} \int_0^\infty l^{k-\lambda-1}e^{-l} \td l -\lim_{l\to0}l^{-\lambda}\pk{Y(l)>x}\\  \label{seri-fini}
 &=& \lambda\sum_{k=1}^\infty \frac{\Gamma(k-\lambda)}{k!}\overline{\MF_\alpha^{*k}}(x)<\infty,
 \EQN
 where the last equality, recalling that $\lambda\in(0,1)$, follows by
 \BQN\label{cov-near0}
 \lim_{l\to0}l^{-\lambda}\pk{Y(l)>x}\leq \lim_{l\to0}l^{-\lambda}(1-e^{-l})=0.
  \EQN
 Next, by a similar argument as used in the proof of Theorem 3.2 in \cite{arendarczyk2012exact}, for any $0<l_0<l_1<\infty$ we have
\BQNY
\pk{L_u^\ast[0,T]>x} &=& \Big(\int_0^{l_0m(u)} +  \int_{l_0m(u)}^{l_1m(u)}  +  \int_{l_1m(u)}^{\infty}\Big)  \pk{L_u^\ast[0,l]>x} \td F_T(l)\\
&=& I_1 +I_2+ I_3,
\EQNY
where
\BQNY
\limsup_{u\to\infty} \frac{I_1}{\pk{T>m(u)}} \leq \limsup_{u\to\infty}
\frac{\int_0^{l_0m(u)}\pk{\sup_{s\in[0,l]}X(s)>u} \td F_T(l)}{\pk{T>m(u)}} \leq \frac{\lambda}{1-\lambda}l_0^{1-\lambda}
\EQNY
and
\BQNY
\limsup_{u\to\infty} \frac{I_3}{\pk{T>m(u)}} \leq \limsup_{u\to\infty}
\frac{\pk{T>l_1m(u)}}{\pk{T>m(u)}} = l_1^{-\lambda}.
\EQNY
Further, in view of Corollary \ref{cor-unc}, for any given $x\geq0$ and arbitrary $\vp>0$, we have the following upper bound
\BQNY
I_2 &=& \int_{l_0}^{l_1} \pk{L_u^*[0,lm(u)]>x} \td F_T(lm(u))\\
&\leq& (1+\vp) \int_{l_0}^{l_1} \pk{Y(l)>x} \td F_T(lm(u))\\
&=& (1+\vp) \Big( \int_{l_0}^{l_1} \pk{T>lm(u)} \td \pk{Y(l)>x} \\
&&\qquad\qquad
- \pk{Y(l_1)>x}\pk{T>l_1m(u)} + \pk{Y(l_0)>x}\pk{T>l_0m(u)} \Big),
\EQNY
which holds for $u$ large enough.
By Potter's Theorem (see, e.g., \cite{Bingham1989}[Theorem 1.5.6]), for sufficiently large $u$ there exists some constant $C>1$ such that
$$\frac{\pk{T>lm(u)}}{\pk{T>m(u)}}\leq C l^{-2\lambda}$$
holds for all $l\in[l_0,l_1]$, and thus by dominated convergence theorem
\BQNY
\limsup_{u\to\infty} \frac{I_2}{\pk{T>m(u)}} \leq (1+\vp)\Big(
\int_{l_0}^{l_1} l^{-\lambda} \td  \pk{Y(l)>x} - \pk{Y(l_1)>x}l_1^{-\lambda} + \pk{Y(l_0)>x}l_0^{-\lambda} \Big).
\EQNY
Similarly, we have the lower bound
\BQNY
\liminf_{u\to\infty} \frac{I_2}{\pk{T>m(u)}} \geq (1-\vp)\Big(
\int_{l_0}^{l_1} l^{-\lambda} \td  \pk{Y(l)>x} - \pk{Y(l_1)>x}l_1^{-\lambda} + \pk{Y(l_0)>x}l_0^{-\lambda} \Big).
\EQNY
Finally, letting $l_0\to0$ and $l_1\to\infty$ in the above bounds,
using \eqref{seri-fini} and \eqref{cov-near0}, and then by the fact that $\vp>0$ was arbitrary, we complete the proof. \QED

\subsection{\prooftheo{theo-sv}}
According to Remark 3.3 in \cite{arendarczyk2012exact}, we know that
\BQNY
\limsup_{u\to\infty}\frac{\pk{L_u^\ast[0,T]>x}}{\pk{T>m(u)}} \leq
\limsup_{u\to\infty}\frac{\pk{\sup_{s\in[0,T]}X(s)>u}}{\pk{T>m(u)}} \leq 1.
\EQNY
Further, by \nekorr{cor-unc}, for arbitrary $l>0$
\BQNY
\pk{L_u^\ast[0,lm(u)]>x} \rightarrow \pk{Y(l)>x}
\EQNY
holds for any $x\geq0$ as $u\to\infty$. Thus, for $T$ with slowly varying tail distribution we get
\BQNY
\liminf_{u\to\infty}\frac{\pk{L_u^\ast[0,T]>x}}{\pk{T>m(u)}} &\geq&
\liminf_{u\to\infty}\frac{\pk{L_u^\ast[0,lm(u)]>x}\pk{T>lm(u)}}{\pk{T>m(u)}}\\
&=& \pk{Y(l)>x},
\EQNY
which converges to $1$ as $l\to\infty$, since by the strong law of large numbers  $Y(l)/l\to\Bal^{-1}(0)>0$.
This completes the proof. \QED

\section{Appendix}\label{appen}
\def\vst{(\vk{s},\vk{t})}
\def\vs{\vk{s}}
\def\vt{\vk{t}}
\def\v0{\vk{0}}
\def\xud{\xi_{u,d}}
\def\Rud{R_{u,d}}
\def\ooint{{\bigcirc}\kern-15pt{\int}\kern-6.5pt{\int}}
\def\chud{\chi_{u,d}}
\def\thud{\theta_{u,d}}
\def\fud{f_{u,d}}
\def\vz{\vk{z}}
\def\Iud{\mathcal{I}_{u,d}}
\def\MI{\mathcal{I}}
\def\vD{\vk{D}}
\def\ve{\vk{e}}

Hereafter, $C_i, i\in \N$ are  positive constants which may be different from line to line.
All vectors are column vectors unless otherwise specified.
As long as it doesn't cause confusion we use $\v0$ to denote  the $2\times1$ column vector or
the $2\times2$ matrix whose entries are all 0's. 
For a given vector (matrix) $\mathbf{Q}$, let $|\mathbf{Q}|$ denote vector (matrix)
with entries equal to absolute value of respective entries of  $\mathbf{Q}$.
\BEL\label{lem-dousum}
Let $X(t),t\geq0$ be a centered stationary Gaussian process with unit variance and covariance function satisfying {\bf A1} and {\bf A3}. If $v(u)$, $\Bal(S,x)$ and $m(u)$ are defined in \eqref{defi-vu}, \eqref{BalSx} and \eqref{def-mu}, respectively,
then for any $A(u)>1$ satisfying
\BQN\label{cond-Au}
\limsup_{u\to\infty}\frac{u^2}{\log A(u)} <\infty
\EQN
we have
\BQNY\label{asy-dousum}
\ \ \ \ \
\lim_{u\to\infty}\sup_{d\geq A(u)}
\abs{\frac{\pk{\sup_{s_1\in[0,S]}X(s_1/v(u))>u, \sup_{s_2\in[0,S]}X(d+s_2/v(u))>u}}{\Psi^2(u)}-\Bal^2(S,0)}=0.
\EQNY
\EEL

\prooflem{lem-dousum} We borrow the argument used in the proof of Theorem 5.1 in \cite{debicki2017approximation}. First, for notational simplicity we define
\BQNY
\xud(\vk{s})=(X(s_1/v(u)),X(d+s_2/v(u)))^T,\quad \vk{s}=(s_1,s_2)\in\vD=[0,S]^2,
\EQNY
and denote by $R_{u,d}(\vk{s},\vk{t})$ the covariance matrix function of $\xud$, i.e.,
\BQNY
\Rud(\vk{s},\vk{t})&=&\Cov(\xud(\vk{s}),\xud(\vk{t}))\\
&=&\E{\xud(\vk{s})\xud(\vk{t})^T}\\
&=&\begin{pmatrix} r(\frac{|t_1-s_1|}{v(u)}) & r(d+\frac{t_2-s_1}{v(u)}) \\ r(d+\frac{s_2-t_1}{v(u)}) & r(\frac{|t_2-s_2|}{v(u)})  \end{pmatrix}
\quad \vk{s},\vk{t}\in\vD.
\EQNY
Then, conditioning on $\xud(\v0)$ we have
 \BQNY
\pk{\exists_{\vk{s}\in\vD} \xi_{u,d}(\vk{s})>\vk{u}}=
\iint_{\R^2}\pk{\exists_{\vk{s}\in\vD} \xi_{u,d}(\vk{s})>\vk{u}| \xud(\v0)=\vk{y}} \phi(y_1,y_2;r(d))\td y_1 \td y_2,
\EQNY
where $\vk{y}=(y_1,y_2)^T$ and $\phi(y_1,y_2;r(d))$ is the density function of
bivariate normal random variable $\xud(\v0)$. By the change of variables $\vk{y}=\vk{u}+\vk{z}/u$ and
using properties of conditional distribution of normal random variable  (see e.g., Chapter 2.2 in \cite{Berman92}), we get
 \BQNY
\pk{\exists_{\vk{s}\in\vD} \xi_{u,d}(\vk{s})>\vk{u}}&=& \frac{e^{-u^2}}{2\pi u^2}
\iint_{\R^2}\pk{\exists_{\vk{s}\in\vD} \chud(\vs) - \thud(\vs,\vz) >\v0 } \fud(\vz) \td z_1 \td z_2\\
&=& \frac{e^{-u^2}}{2\pi u^2} \iint_{\R^2} \Iud(\vz)\fud(\vz) \td z_1 \td z_2,
\EQNY
where
$
\Iud(\vz):=\pk{\exists_{\vk{s}\in\vD} \chud(\vs) - \thud(\vs,\vz) >\v0 }
$
with
$$ \chud(\vs) = u\LT( \xud(\vs)-\Rud(\vs,\v0)\Rud^{-1}(\v0,\v0)\xud(\v0)\RT),\quad \vs\in\vD,$$
$$\thud(\vs,\vz)=u^2\LT(\vk{1}-\Rud(\vs,\v0)\Rud^{-1}(\v0,\v0)(\vk{1}+\vz/u^2)\RT),\quad \vs\in\vD, \vk{z}\in\R^2$$
and
$$\fud(\vk{z})=\frac{1}{\sqrt{1-r^2(d)}}
\exp\LT(\frac{1}{1+r(d)}(u^2r(d)-z_1-z_2)-\frac{(z_1-r(d)z_2)^2}{2u^2(1-r^2(d))}-\frac{z_2^2}{2u^2}\RT),\quad \vk{z}\in\R^2.$$
Consequently, in order to show the claim it suffices to prove that
for
$W_\alpha(\vs)=(\sqrt{2}B_\alpha^{(1)}(s_1)-s_1^\alpha, \sqrt{2}B_\alpha^{(2)}(s_2)-s_2^\alpha)^T$,
$\vs\in\vD$,
where $B_\alpha^{(1)}$ and $B_\alpha^{(2)}$ are two independent fBm's with Hurst index $\alpha/2$,
\BQN \nonumber
\lefteqn{\lim_{u\to\infty}\sup_{d\geq A(u)}
\abs{\iint_{\R^2}\Iud(\vz) \fud(\vz) \td z_1 \td z_2-\Bal^2(S,0)}}\\ \nonumber
&&=\lim_{u\to\infty}\sup_{d\geq A(u)}
\abs{\iint_{\R^2} \bigg(\Iud(\vz) \fud(\vz)  -
\pk{\exists_{\vs\in\vD} W_\alpha(\vs)+\vz>\v0}  e^{-z_1-z_2}\bigg) \td z_1 \td z_2 }\\ \label{asy-dousum-ma}
&&=\lim_{u\to\infty}\sup_{d\geq A(u)}
\abs{\iint_{\R^2} \big(\Iud(\vz) \fud(\vz)  - \MI(\vz) e^{-z_1-z_2}\big) \td z_1 \td z_2 }=0,
\EQN
with
\begin{eqnarray}
\MI(\vz)&:=&\pk{\exists_{\vs\in\vD} W_\alpha(\vs)+\vz>\v0} \label{prod0}\\
&=&
\pk{\sup_{s_\in[0,S]} \sqrt{2}B_\alpha(s)-s^\alpha >-z_1}
\pk{\sup_{s\in[0,S]} \sqrt{2}B_\alpha(s)-s^\alpha >-z_2}.\label{prod}
\end{eqnarray}

\quad For any $\vs,\vt\in\vD$,
\BQN \nonumber
\lefteqn{\Cov(\chud(\vs),\chud(\vt))}\\ \nonumber
&&=u^2
\Cov\big(\xud(\vs)-\Rud(\vs,\v0)\Rud^{-1}(\v0,\v0)\xud(\v0),\,\, \xud(\vt)-\Rud(\vt,\v0)\Rud^{-1}(\v0,\v0)\xud(\v0)\big)\\ \nonumber
&&=u^2 \{ \Rud(\vs,\vt)- \Rud(\vs,\v0)\Rud^{-1}(\v0,\v0)\Rud(\v0,\vt) \}\\ \label{covchi}
&&=u^2 \{ (\Rud(\vs,\vt)-E) + (E-\Rud(\vs,\v0)) + \Rud(\vs,\v0)(E-\Rud^{-1}(\v0,\v0)\Rud(\v0,\vt)) \}£¬
\EQN
where $E$ is the $2\times2$ identity matrix.

Since $A(u)>1$ satisfying \eqref{cond-Au} tends to $\infty$ as $u\to\infty$, then by {\bf A3} we have
\BQN\label{rcovas1}
\lim_{u\to\infty }\sup_{d>A(u),\vs\in\vD} \abs{ \Rud(\vs,\v0)-E}=\v0
\EQN
and
\BQN\label{u2rlimit}
\lim_{u\to\infty } u^2\sup_{d\geq A(u)}\abs{r(d)}\leq \lim_{u\to\infty }\frac{u^2}{\log A(u)} \sup_{d\geq A(u)}\abs{r(d)}\log d =0.
\EQN
Therefore,
 \BQN\label{rcovas2}
 \lim_{u\to\infty }\sup_{d\geq A(u)} u^2(E-\Rud^{-1}(\v0,\v0))
 =\lim_{u\to\infty }\sup_{d\geq A(u)}
 \frac{u^2r(d)}{1-r^2(d)} \begin{pmatrix} -r(d) & 1 \\ 1 & -r(d) \end{pmatrix} =\v0.
 \EQN
  Note that by {\bf A1}, \eqref{defi-vu} and the Uniform Convergence Theorem (see, e.g., Theorem 1.5.2 in \cite{Bingham1989}) we get
\BQNY
&&\quad\lim_{u\to\infty }\sup_{s\in[0,S]} \abs{u^2(1-r(s/v(u))) -s^\alpha}\\
&&\leq\lim_{u\to\infty }\sup_{s\in[0,S]} \abs{\frac{1-r(s/v(u))}{1-r(1/v(u))}-s^\alpha }
+ \lim_{u\to\infty }\sup_{s\in[0,S]} \abs{\frac{1-r(s/v(u))}{1-r(1/v(u))}}
\abs{u^2(1-r(1/v(u)))-1}\\
&&=0.
\EQNY
 Consequently,
\BQN\label{rcovas3}
\lim_{u\to\infty }\sup_{d\geq A(u),\vs,\vt\in\vD} \abs{u^2(E-\Rud(\vs,\vt))
- \begin{pmatrix} \abs{t_1-s_1}^\alpha & 0 \\ 0 & \abs{t_2-s_2}^\alpha \end{pmatrix}} =\v0.
\EQN
Similarly,
\BQN\nonumber \label{rcovas4}
&&\quad\lim_{u\to\infty }\sup_{d\geq A(u),\vs,\vt\in\vD}
\abs{u^2 \Rud(\vs,\v0)\big(E-\Rud^{-1}(\v0,\v0)\Rud(\v0,\vt)\big)
- \begin{pmatrix} \abs{t_1}^\alpha & 0 \\ 0 & \abs{t_2}^\alpha \end{pmatrix}}\\ \nonumber
&&\leq \lim_{u\to\infty }\sup_{d\geq A(u),\vs,\vt\in\vD} \abs{u^2(E-\Rud(\v0,\vt))
- \begin{pmatrix} \abs{t_1}^\alpha & 0 \\ 0 & \abs{t_2}^\alpha \end{pmatrix}} \\ \nonumber
&&\quad + \lim_{u\to\infty }\sup_{d\geq A(u),\vs,\vt\in\vD}
\abs{( \Rud(\vs,\v0) -E) [u^2(E-\Rud(\v0,\vt))] }\\ \nonumber
&&\quad + \lim_{u\to\infty }\sup_{d\geq A(u),\vs,\vt\in\vD}
\abs{ \Rud(\vs,\v0) [u^2(E-\Rud^{-1}(\v0,\v0))] \Rud(\v0,\vt) }\\
&&=\v0,
\EQN
where in the last equality we have used \eqref{rcovas1}, \eqref{rcovas2} and \eqref{rcovas3}.

Substituting \eqref{rcovas3} and \eqref{rcovas4} into \eqref{covchi} gives
\BQN
 {\small \lim_{u\to\infty }\sup_{d\geq A(u),\vs,\vt\in\vD}
\Bigg|\Cov(\chud(\vs),\chud(\vt))- 
\begin{pmatrix} |s_1|^\alpha+|t_1|^\alpha-\abs{t_1-s_1}^\alpha & 0 \\ 0 & |s_2|^\alpha+|t_2|^\alpha-\abs{t_2-s_2}^\alpha \end{pmatrix}\Bigg|} &&
\nonumber\\
&=&\v0.
\label{weak-cov-chi}
\EQN
 Hence, the finite-dimensional distributions of $\chud$ converge to that
 of $\{\sqrt{2}B_\alpha(\vs), \vs\in\vD\}$ uniformly with respect to $d\geq A(u)$, where
 $B_\alpha(\vs)=(B_\alpha^{(1)}(s_1),B_\alpha^{(2)}(s_2))^T.$

\quad Let $C(\vD)$ denote the Banach space of all continuous functions on $\vD$ equipped with sup-norm, we now show that the measures on $C(\vk{D})$ induced by $\{\chud(\vs),\vs\in\vk{D},d\geq A(u)\}$ are uniformly tight for large $u$.
In fact, since $\E{\xud(\vs)|\xud(\v0)}=\Rud(\vs,\v0)\Rud^{-1}(\v0,\v0)\xud(\v0)$ then
\BQNY
&&\quad\E{(\xud(\vs)-\xud(\vt))^T\LT(\Rud(\vs,\v0)-\Rud(\vt,\v0)\RT)\Rud^{-1}(\v0,\v0)\xud(\v0)}\\
&&=\E{\E{(\xud(\vs)-\xud(\vt))^T|\xud(\v0)}\LT(\Rud(\vs,\v0)-\Rud(\vt,\v0)\RT)\Rud^{-1}(\v0,\v0)\xud(\v0)}\\
&&=\E{\xud^T(\v0)\Rud^{-1}(\v0,\v0)(\Rud(\v0,\vs)-\Rud(\v0,\vt))
(\Rud(\vs,\v0)-\Rud(\vt,\v0))\Rud^{-1}(\v0,\v0)\xud(\v0)}
\EQNY
and thus
\BQNY
&&\quad\E{\norm{\chud(\vs)-\chud(\vt)}^2}\\
&&=u^2\Big[\E{\norm{\xud(\vs)-\xud(\vt)}^2} -2\E{(\xud(\vs)-\xud(\vt))^T(\Rud(\vs,\v0)-\Rud(\vt,\v0))\Rud^{-1}(\v0,\v0)\xud(\v0)}\\
&&\quad + \E{\xud^T(\v0)\Rud^{-1}(\v0,\v0)(\Rud(\v0,\vs)-\Rud(\v0,\vt))
(\Rud(\vs,\v0)-\Rud(\vt,\v0))\Rud^{-1}(\v0,\v0)\xud(\v0)} \Big]\\
&&=  u^2\Big[\E{\norm{\xud(\vs)-\xud(\vt)}^2}
 - \E{\norm{(\Rud(\vs,\v0)-\Rud(\vt,\v0))\Rud^{-1}(\v0,\v0)\xud(\v0)}^2} \Big]\\
&& \leq u^2 E{\norm{\xud(\vs)-\xud(\vt)}^2} \\
&& =2u^2[1-r(\abs{t_1-s_1}/v(u)) + 1-r(\abs{t_2-s_2}/v(u)) ].
\EQNY
Moreover, by \eqref{defi-vu} and Potter's Theorem (see, e.g., \cite{Bingham1989}[Theorem 1.5.6]), for large enough $u$
there exists some constant $C>1$ such that
\BQNY\label{Potter-bound}
u^2(1-r(s/v(u)))=u^2 (1-r(1/v(u)))\frac{1-r(s/v(u))}{1-r(1/v(u))}\leq C\abs{s}^{\alpha/2}
\EQNY
holds for all $s\in[0,S]$. Hence, for large enough $u$ we get
\BQN\label{incvarich}
\sup_{d\geq A(u)}\E{\norm{\chud(\vs)-\chud(\vt)}^2} \leq 2C(|t_1-s_1|^{\alpha/2}+|t_2-s_2|^{\alpha/2})
\EQN
for any $\vs,\vt\in\vD$, implying the uniform tightness of the measures induced by $\{\chud(\vs),\vs\in\vk{D},d\geq A(u)\}$.
This together with \eqref{weak-cov-chi} implies that $\chud$ converges weakly, as $u\to\infty$, to $\sqrt{2}B_\alpha(\vs),\vs\in\vD$ uniformly for $d\geq A(u)$.

Further, by \eqref{rcovas2} and \eqref{rcovas3} we have
\BQN\label{lim-thud}
\lefteqn{\lim_{u\to\infty}\sup_{d\geq A(u),\vs\in\vD}
\abs{u^2\big(E-\Rud(\vs,\v0)\Rud^{-1}(\v0,\v0)\big)
-\begin{pmatrix} \abs{s_1}^\alpha & 0 \\ 0 & \abs{s_2}^\alpha \end{pmatrix}}}\\ \nonumber
&&\leq \lim_{u\to\infty}\sup_{d\geq A(u),\vs\in\vD} \abs{u^2(E-\Rud(\vs,\v0))
-\begin{pmatrix} \abs{s_1}^\alpha & 0 \\ 0 & \abs{s_2}^\alpha \end{pmatrix}}\\ \nonumber
&&\quad + \lim_{u\to\infty}\sup_{d\geq A(u),\vs\in\vD} \abs{\Rud(\vs,\v0)[u^2(E-\Rud^{-1}(\v0,\v0))]}\\ \nonumber
&&=\v0
\EQN
and thus for any $\vz\in\R^2$
\BQN\nonumber
\lefteqn{\lim_{u\to\infty}\sup_{d\geq A(u),\vs\in\vD}\abs{\thud(\vs)-(|s_1|^\alpha-z_1,|s_2|^\alpha-z_2)^T}}\\ \nonumber
&&\leq \lim_{u\to\infty}\sup_{d\geq A(u),\vs\in\vD}
\abs{\Bigg[u^2\big(E-\Rud(\vs,\v0)\Rud^{-1}(\v0,\v0)\big) -
\begin{pmatrix} \abs{s_1}^\alpha & 0 \\ 0 & \abs{s_2}^\alpha \end{pmatrix}\Bigg]
\begin{bmatrix} 1 \\ 1 \end{bmatrix}}\\ \nonumber
&&\quad + \quad \lim_{u\to\infty}\sup_{d\geq A(u),\vs\in\vD}
\abs{\big(E-\Rud(\vs,\v0)\Rud^{-1}(\v0,\v0)\big)\vz}\\ \nonumber
&&=\v0.
\EQN
Therefore, for each $\vz\in\R^2$, the probability measures on $C(\vk{D})$ induced by $\{\chud(\vs)-\thud(\vs,\vz),\vs\in\vD\}$
converge weakly, as $u\to\IF$, to that induced by $\{W_\alpha(\vs)+ \vz,\vs\in\vD\}$ uniformly with respect to $d\geq A(u)$.
Then, by the continuous mapping theorem, (\ref{prod}) and the fact that
the set of discontinuity points of cumulative distribution function of
$\sup_{s_\in[0,S]} \sqrt{2}B_\alpha(s)-s^\alpha$
consists of at most of one point
(see, e.g., Theorem 7.1 in \cite{AzW09} or related Lemma 4.4 in \cite{DHW19}), we get
\BQNY
\lim_{u\to\infty}\sup_{d\geq A(u)}\abs{\Iud(\vz)-\MI(\vz)}=0
\EQNY
for almost all $\vz\in \R^2$, where $\MI(\vz)$ is defined in \eqref{prod0}. Further, by \eqref{u2rlimit} we know
$$\lim_{u\to\infty}\sup_{d\geq A(u)}\abs{f_{u,d}(\vz)-e^{-z_1-z_2}}=0,\quad\forall\, \vz\in\R^2,$$
and thus for almost all $\vz\in \R^2$
\BQN\label{weakcov-I}
\lim_{u\to\infty}\sup_{d\geq A(u)}\abs{\Iud(\vz)f_{u,d}(\vz)-\MI(\vz)e^{-z_1-z_2}}=0.
\EQN
\quad Therefore, to verify \eqref{asy-dousum-ma}, we have to put the limit into integral. In the following, we
 look for an integrable upper bound for $\sup_{d\geq A(u)}\Iud(\vz)f_{u,d}(\vz)$. We first give a lower bound for
$\inf_{d\geq A(u),\vs\in\vD}\thud(\vs,\vz)$. Let $\vp(<1/2)$ be a positive constant.
In view of \eqref{lim-thud}, we know that, for sufficiently large $u$
\BQNY
\sup_{d\geq A(u),\vs\in\vD}\abs{E-\Rud(\vs,\v0)\Rud^{-1}(\v0,\v0)}\leq
\begin{pmatrix} \vp & \vp \\ \vp & \vp\end{pmatrix},
\EQNY
and thus
\begin{eqnarray*}
\inf_{d\geq A(u),\vs\in\vD}\thud(\vs,\vz)
&=&  \inf_{d\geq A(u),\vs\in\vD}
\{u^2(E-\Rud(\vs,\v0)\Rud^{-1}(\v0,\v0))\vk{1}
+(E-\Rud(\vs,\v0)\Rud^{-1}(\v0,\v0))\vz-\vz\}\\
&\geq&
 -\vk{1}-\vz + \inf_{d\geq A(u),\vs\in\vD} \{(E-\Rud(\vs,\v0)\Rud^{-1}(\v0,\v0))\vz\} \\
& \geq& -\vk{1}-\vz - \begin{pmatrix} \vp & \vp \\ \vp & \vp\end{pmatrix} \abs{\vz}
:=h(\vz),\quad \vz\in\R^2. \label{solution-W}
\end{eqnarray*}
Let $\{\ve_k,k=1,2,3\}$ denotes $(1,1)^T$, $(0,1)^T$ and $(1,0)^T$, respectively. By Cauchy-Schwartz inequality and \eqref{incvarich}, for large enough $u$
\BQN\nonumber
\sup_{d\geq A(u)}\E{ \LT(\ve_k^T(\chud(\vs)-\chud(\vt))\RT)^2 } &\leq& \sup_{d\geq A(u)}
2 \E{\norm{\chud(\vs)-\chud(\vt)}^2} \\ \label{upp-var-chud}
&\leq& 4C(|t_1-s_1|^{\alpha/2}+|t_2-s_2|^{\alpha/2}),\quad k=1,2,3
\EQN
holds for any $\vs,\vt\in\vD$. Thus, by Sudakov-Fernique inequality (see, e.g., \cite{adler1990introduction}[Theorem 2.9]),
 we have
\BQN\label{S-F-ine}
\sup_{d\geq A(u)}\E{\sup_{\vk{s}\in \vk{D}} \ve_k^T\chud(\vs) } \leq
\E{\sup_{\vk{s}\in \vk{D}} \sum_{i=1}^2 2\sqrt{C}B_{\alpha/2}^{(i)}(s_i) }:=C_1<\infty,\quad k=1,2,3,
\EQN
where $B_{\alpha/2}^{(i)}$'s are independent fBm's with Hurst index $\alpha/4$.
Then, for all large enough $u$,
\BQN\nonumber
\sup_{d\geq A(u)} \Iud(\vz)&=& \sup_{d\geq A(u)}\pk{\exists_{\vk{s}\in\vD} \chud(\vs) - \thud(\vs,\vz) >\v0 }\\ \nonumber
&\leq& \sup_{d\geq A(u)} \pk{\exists_{\vs\in\vD} \chud(\vs) > \inf_{d\geq A(u),\vs\in\vD} \thud(\vs,\vz) }\\ \nonumber
&\leq& \sup_{d\geq A(u)} \pk{\sup_{\vs\in\vD} \ve_k^T\chud(\vs) > \ve_k^T h(\vz) }\\ \label{BIS-xud}
&\leq& \sup_{d\geq A(u)}
\exp\LT( - \frac{\LT(\ve_k^T h(\vz)-\E{\sup_{\vk{s}\in \vk{D}} \ve_k^T\chud(\vs) }\RT)^2} {2\Var_{\vs\in\vD} \ve_k^T\chud(\vs)} \RT)\\ \nonumber
&\leq&
\exp\LT( - C_2\LT(\ve_k^T h(\vz)-C_1\RT)^2 \RT),\quad \vz\in \vk{Z}_k, k=1,2,3,
\EQN
where \eqref{BIS-xud} follows from Borell-TIS inequality (see, e.g., \cite{AdlerTaylor}[Theorem 2.1.1]), the last inequality follows by \eqref{upp-var-chud}-\eqref{S-F-ine} with $C_2=(16CS^{\alpha/2})^{-1}$, and
$$\vk{Z}_1=\{(z_1,z_2)| z_1<0,z_2<0,(2\vp-1)(z_1+z_2)>2+C_1\},$$
$$\vk{Z}_2=\{(z_1,z_2)| z_1>0,z_2<0,(\vp-1)z_2-\vp z_1>1+C_1\},$$
$$\vk{Z}_3=\{(z_1,z_2)| z_1<0,z_2>0,(\vp-1)z_1-\vp z_2>1+C_1\}.$$
 Therefore,
 \BQNY
 \sup_{d\geq A(u)} \Iud(\vz) \leq g(\vz) :=
 \left\{\begin{array}{ll}
\exp\LT( - C_2\LT(\ve_k^T h(\vz)-C_1\RT)^2 \RT),\quad& \vz\in \vk{Z}_k, k=1,2,3,\\
1,\quad& \quad \vz\in \R^2\backslash\bigcup_{k=1}^3\vk{Z}_k,\\
\end{array}\right.\\
 \EQNY
 holds for sufficiently large $u$. Moreover,  by \eqref{u2rlimit}
  \BQNY
&&\quad \sup_{d\geq A(u)} f_{u,d}(\vz)e^{z_1+z_2}\\
&& = \sup_{d\geq A(u)} \frac{1}{\sqrt{1-r^2(d)}}\exp\LT(\frac{u^2r(d)}{1+r(d)}+
  \frac{ 2u^2r(d)(1-r(d))(z_1+z_2) -(z_1^2-2r(d)z_1z_2+z_2^2)}{2u^2(1-r^2(d))}\RT)\\
 && \leq \frac{3}{2}\sup_{d\geq A(u)} \exp\LT(\frac{u^2r(d)}{1+r(d)}+
  \frac{ -\frac{1+r(d)}{2}(z_2-z_1)^2 -\frac{1-r(d)}{2}(z_1+z_2-2u^2r(d))^2 +2u^4r^2(d)(1-r(d)) }{2u^2(1-r^2(d))}\RT)\\
 && \leq \frac{3}{2}\sup_{d\geq A(u)} e^{u^2r(d)} \leq 2,\quad \vz\in\R^2
 \EQNY
 holds for all large enough $u$, and thus
  \BQNY
  \sup_{d\geq A(u)}\Iud(\vz)f_{u,d}(\vz) \leq 2g(\vz)e^{-z_1-z_2},\quad \vz\in\R^2.
  \EQNY
   We now show that
  $g(\vz)e^{-z_1-z_2}$ is integrable on $\R^2$.
   In fact,
 \BQNY
 \iint_{\R^2} g(\vz)e^{-z_1-z_2} \td z_1 \td z_2 =\LT( \iint_{\vk{Z}_1} +\iint_{\vk{Z}_2} +\iint_{\vk{Z}_3}
 +\iint_{\R^2\backslash\bigcup_{k=1}^3\vk{Z}_k}\RT) g(\vz)e^{-z_1-z_2} \td z_1 \td z_2,
  \EQNY
  where
  \BQNY
\iint_{\vk{Z}_1} g(\vz)e^{-z_1-z_2} \td z_1 \td z_2 &\leq&
\int_{-\infty}^0 \int_{-\infty}^0
\exp\LT(-C_2\LT( (2\vp-1)z_1+(2\vp-1)z_2-2-C_1 \RT)^2-z_1-z_2\RT)
 \td z_1 \td z_2\\
 &\leq& \LT(\int_{-\infty}^0 \exp\LT( -C_2(2\vp-1)^2z_1^2+\LT(2C_2(2\vp-1)(2+C_1)-1\RT)z_1 \RT)  \td z_1\RT)^2\\
 &<& \infty,
  \EQNY
   \BQNY
\iint_{\vk{Z}_2} g(\vz)e^{-z_1-z_2} \td z_1 \td z_2 &= &
  \iint_{\vk{Z}_3} g(\vz)e^{-z_1-z_2} \td z_1 \td z_2\\
  &=& \int^{\infty}_0 \LT(\int_{-\infty}^{\frac{\vp z_1+1+C_1}{(\vp-1)}}
\exp\LT(-C_2\LT( (\vp-1)z_2-\vp z_1-1-C_1 \RT)^2-z_2\RT)
 \td z_2\RT) e^{-z_1} \td z_1\\
 &=& \frac{e^{\frac{1+C_1}{1-\vp}}}{1-\vp}\int^{\infty}_0  e^{(\frac{\vp}{1-\vp}-1)z_1} \td z_1
 \int_{0}^{\infty}
\exp\LT(-C_2z_2^2-\frac{z_2}{\vp-1}\RT) \td z_2
 < \infty,
  \EQNY
  since $\vp<1/2$, and
\BQNY
\iint_{\R^2\backslash\bigcup_{k=1}^3\vk{Z}_k} g(\vz)e^{-z_1-z_2} \td z_1 \td z_2
&\leq& \LT(\iint_{z_1<0,z_2<0, z_1+z_2\geq\frac{2+C_1}{2\vp-1}}
+ 2\int_0^\infty\int^\infty_{\frac{\vp z_1+C_1+1}{\vp-1}} \RT)
e^{-z_1-z_2} \td z_1 \td z_2\\
&\leq& \LT(\frac{2+C_1}{1-2\vp}\RT)^2e^{\frac{2+C_1}{1-2\vp}} +
2e^{\frac{1+C_1}{1-\vp}}\int^{\infty}_0  e^{(\frac{\vp}{1-\vp}-1)z_1} \td z_1
 < \infty.
\EQNY
\quad Consequently, \eqref{asy-dousum-ma} follows by the dominated convergence theorem and \eqref{weakcov-I}. This completes the proof. \QED

\BEL\label{thasysmu}
Let $X(t),t\geq0$ be a centered stationary Gaussian process with unit variance and covariance function satisfying
{\bf A1} and {\bf A3}. Let $v(u)$, $\Bal(x)$ and $m(u)$ be defined in \eqref{defi-vu}, \eqref{def-B-H-x} and \eqref{def-mu} respectively.
Then for $A(u)$ such that 
\BQN\label{vAm}
\lim_{u\to\infty}A(u)v(u)=\IF\quad  \textrm{and}\quad \lim_{u\to\infty}\frac{A(u)}{m(u)}=0
\EQN
and any $x\geq0$ we have
\BQN\label{asysojsmu}
\pk{L_u^\ast[0,A(u)]>x} \sim \MB_\alpha(x)A(u)v(u)\Psi(u),\quad u\to\infty.
\EQN
\EEL

\prooflem{thasysmu} We follow the argument used in the proof of Theorem 2.1 in \cite{debicki2017approximation}.
Let $A(u)$ satisfy \eqref{vAm}, for any $S>1$ define
$$\Delta_k=[k S/v(u),(k+1) S/v(u)],\quad  k=0 ,\ldots,N_u
$$
with $N_u=\lfloor A(u)v(u) /S \rfloor$, i.e., the integer part of $A(u)v(u) /S $. By stationarity of $X$, we have for all $u$ positive and $x\ge 0$
\BQNY\label{low-upp-I1}
 I_1(u) \leq \pk{ L_{u}^*[0,T] >x} \leq  I_2(u),
\EQNY
 where
\BQNY
&& I_1(u) =(N_u-1)\pk{ L_{u}^*\Delta_0>x  } - \sum_{0\leq i<k\leq N_u-1}
q_{i,k}(u),\\
&&  I_2(u) = (N_u+1)\pk{ L_{u}^*\Delta_0>x } + \sum_{0\leq i<k\leq N_u}
q_{i,k}(u),
\EQNY
with
$q_{i,k}(u)=\pk{ \sup_{t\in\Delta_i} X(t)>u, \sup_{t\in\Delta_k} X(t)>u }.
$
By Theorem 5.1 in \cite{debicki2017approximation} and \eqref{def-B-H-x}, we have
\BQN\label{sojsum}
\lim_{S\to\IF}\lim_{u\to\IF} \frac{N_u\pk{ L_{u}^*\Delta_0>x }}{A(u)v(u)\Psi(u)} =\MB_\alpha(x)
\EQN
 for any $x\geq0$. Therefore, it suffices to show that the double sum
is negligible with respect to $A(u)v(u)\Psi(u)$ as $u\to\IF$ and then as $S\to\IF$.

\quad Let $\vp^\ast(<2)$ be the positive root of equation $x^2-(2-\alpha)x-\frac{3}{2}\alpha=0$ and put
$\beta=\inf_{t\geq1}\{1-r(t)\},$
which by {\bf A3} is positive. Define
\BQNY
A_0(u)=0,\,A_1(u)=u^{\frac{\vp^*-2}{\alpha}}\wedge A(u),\, A_2(u)=1\wedge A(u),\,
A_3(u)=e^{\beta u^2/8}\wedge A(u),\,A_4(u)=A(u)
\EQNY
 and
$$\Lambda_l(u)=\{(i,k):1\leq i+1<k\leq N_u, A_{l-1}(u)<(k-i-1)S/v(u)\leq A_l(u),\,l=1,2,3,4\}.$$
Then
\BQN\label{dou-sum-ls}
\sum_{0\leq i<k\leq N_u} q_{i,k}(u) &=&
\sum_{0\leq i< N_u} q_{i,i+1}(u) + \sum_{l=1}^4\sum_{(i,k)\in \Lambda_l(u)} q_{i,k}(u) \\ \nonumber
&:=&Q_0(u)+\sum_{l=1}^4 Q_l(u).
\EQN
According to (4.7)-(4.9) in \cite{debicki2017approximation} we know that
 \BQN\label{limQ2}
 \limsup_{u\to\IF}\frac{Q_2(u)}{A(u)v(u)\Psi(u)}=0,
 \EQN
  \BQN\label{limQ1}
 \lim_{S\to\IF}\limsup_{u\to\IF}\frac{Q_1(u)}{A(u)v(u)\Psi(u)}=0,
 \EQN
 and
   \BQN\label{limQ0}
 \lim_{S\to\IF}\limsup_{u\to\IF}\frac{Q_0(u)}{A(u)v(u)\Psi(u)}=0.
 \EQN
 Next, by stationarity of $X$, for sufficiently large $u$
\BQNY
\sup_{(i,k)\in \Lambda_3(u)}\E{\sup_{s\in\Delta_i,t\in\Delta_k}(X(s)+X(t))}\leq
2\E{\sup_{s\in[0,1]}X(s)}=:C_3<\IF,
\EQNY
\BQNY
\sup_{(i,k)\in \Lambda_3(u), s\in\Delta_i,t\in\Delta_k} \Var(X(s)+X(t))&=&4-2\inf_{(i,k)\in \Lambda_3(u),(s,t)\in\Delta_i\times\Delta_k}\{1-r(t-s)\}\\
&\leq& 4-2\beta.
\EQNY
Then, by Borell-TIS inequality we have for large enough $u$
\BQNY\nonumber
\sup_{(i,k)\in \Lambda_3(u)}q_{i,k}(u)&\leq&\sup_{(i,k)\in \Lambda_3(u)} \pk{ \sup_{s\in\Delta_i,t\in\Delta_k} X(s)+X(t)>2u}\\
&\leq& \exp\LT(-\frac{(2u-C_3)^2}{2(4-2\beta)}\RT)\\
&\leq& \exp\LT(-\frac{1+\beta/2}{2}(u-C_3/2)^2\RT),
\EQNY
and thus
\BQN\nonumber
\limsup_{u\to\IF}\frac{Q_3(u)}{A(u)v(u)\Psi(u)}
& \leq&
\limsup_{u\to\IF}\frac{N_u A_3(u)v(u)}{SA(u)v(u)\Psi(u)}\exp\LT(-\frac{1+\beta/2}{2}(u-C_3/2)^2\RT) \\ 
& \leq& \limsup_{u\to\IF} \frac{\sqrt{2\pi}u v(u)}{S^2} \exp\LT( -\frac{\beta}{8}u^2+C_3u(1+\beta/2)\RT) \label{limQ3}
 =0.
\EQN
Further, since $e^{\beta u^2/8}$ satisfies \eqref{cond-Au}, then by \nelem{lem-dousum} and stationarity of $X$,
$$Q_4(u)\leq 2N^2_u\Psi^2(u)\Bal^2(S,0)$$
holds for $u$ sufficiently large. 
Therefore,
\BQN\nonumber
\limsup_{u\to\IF}\frac{Q_4(u)}{A(u)v(u)\Psi(u)}
 &\leq&
\limsup_{u\to\IF}\frac{2N_u^2\Psi^2(u)\Bal^2(S,0) }{A(u)v(u)\Psi(u)}\\ \nonumber
& \leq& \limsup_{u\to\IF} \frac{2\MB_{\alpha}^2(S,0)}{S^2\Bal(0)}\frac{A(u)}{m(u)} \label{limQ4}
 =0,
\EQN
where the last equality follows by \eqref{vAm}.

\quad Consequently, substituting  \eqref{limQ2}-\eqref{limQ4} into \eqref{dou-sum-ls} yields
\BQNY\label{ds-neg-ls}
\lim_{S\to\IF}\limsup_{u\to\IF}\frac1{A(u)v(u)\Psi(u)}\sum_{0\leq i<k\leq N_u}
q_{i,k}(u) = 0,
\EQNY
which together with \eqref{sojsum} completes the proof. \QED

\BK\label{corsm}
If $X$, $v(u)$, $\Bal(x)$, $m(u)$ and $A(u)$ are given as in \nelem{thasysmu}, then for any $x\geq0$ and $\vp\in(0,1)$
there exists $\delta>0$ such that
\BQN\label{asysojinf}
\liminf_{u\to\IF}\inf_{t\in[A(u),\delta m(u)]}\frac{\pk{L_u^*[0,t]>x}}{t\MB_\alpha(x)v(u)\Psi(u)} \geq 1-\vp
\EQN
and
\BQN\label{asysojsup}
\limsup_{u\to\IF}\sup_{t\in[A(u),\delta m(u)]}\frac{\pk{L_u^*[0,t]>x}}{t\MB_\alpha(x)v(u)\Psi(u)} \leq 1+\vp.
\EQN
\EK

\proofkorr{corsm} Let $x\geq0$ be fixed, recalling \eqref{def-cp} we have that,
 for arbitrary $\vp>0$, there exists some $\delta>0$ such that
\BQN\label{vpdelta}
(1-\vp/4) \leq \frac{\pk{Y(t)>x}}{t \overline{\MF_\alpha}(x)} \leq (1+\vp/4),\quad t\in(0,\delta).
\EQN
For such $\vp$ and $\delta$, suppose that \eqref{asysojinf} does not hold. Then, there exist two sequences $\{u_n,n\in \N\}$ and $\{t_n,n\in\N\}$ such that $u_n\to\IF$ as $n\to\IF$ and
\BQN\label{asysojinfcotr}
\frac{\pk{L_{u_n}^*[0,t_n]>x}}{t_n\MB_\alpha(x)v(u_n)\Psi(u_n)} < 1-\vp,\quad t_n\in[A(u_n),\delta m(u_n)], n\in \N.
\EQN
Putting ${\hat t}_n=t_n/m(u_n)$, by \eqref{def-mu} and \eqref{def-hf-alp}, we get 
\begin{equation}
 \tag{\ref{asysojinfcotr}$'$} \label{asysojinfcotrv2}
\frac{\pk{L_{u_n}^*[0,{\hat t}_nm(u_n)]>x}}{{\hat t}_n\overline{\MF_\alpha}(x)} < 1-\vp,\quad {\hat t}_n\in[A(u_n)/m(u_n),\delta], n\in \N.
\end{equation}
Since sequence $\{{\hat t}_n,n\in\N\}$ is bounded, then there exists a convergent subsequence
$\{{\hat t}_{n_k},k\in\N\}$ such that 
$\lim_{k\to\IF}{\hat t}_{n_k}\ge0$. 
If $\lim_{k\to\IF}{\hat t}_{n_k}>0$, then by \nekorr{cor-unc}
\BQNY
\frac{\pk{L_{u_{n_k}}^*[0,{\hat t}_{n_k}m(u_{n_k})]>x}}{\pk{Y({\hat t}_{n_k})>x}} > 1-\vp/4
\EQNY
holds for sufficiently large $k$, which together with \eqref{vpdelta} implies
\BQNY
\frac{\pk{L_{u_{n_k}}^*[0,{\hat t}_{n_k}m(u_{n_k})]>x}}{{\hat t}_{n_k}\overline{\MF_\alpha}(x)} > (1-\vp/4)^2.
\EQNY
This however contradicts 
\eqref{asysojinfcotrv2}.
If $\lim_{k\to\IF}{\hat t}_{n_k}=0$, then 
\BQNY
\lim_{k\to\infty}t_{n_k}v(u_{n_k})\geq \lim_{k\to\infty}A(u_{n_k})v(u_{n_k})=\IF\quad  \textrm{and}\quad \lim_{k\to\infty}\frac{t_{n_k}}{m(u_{n_k})}=\lim_{k\to\IF}{\hat t}_{n_k}=0,
\EQNY
and thus by \nelem{thasysmu}
\BQNY
\frac{\pk{L_{u_{n_k}}^*[0,t_{n_k}]>x}}{t_{n_k}\MB_\alpha(x)v(u_{n_k})\Psi(u_{n_k})} > 1-\vp/4
\EQNY
holds for sufficiently large $k$. This contradicts \eqref{asysojinfcotr}.
An analogous argument can be used to verify \eqref{asysojsup}. This completes the proof. \QED
\\

{\bf Acknowledgement}:
The authors would like to thank Enkelejd Hashorva for his numerous valuable remarks
on all the steps of preparation of the manuscript.
K.D. was partially supported by NCN Grant No  2018/31/B/ST1/00370 (2019-2022).
X.P. thanks National Natural Science Foundation of China (11701070,71871046) for partial financial support.
Financial support from the Swiss National Science Foundation Grant 200021-175752/1 is also kindly acknowledged.

\bibliographystyle{ieeetr}
\bibliography{BermanC}
\end{document}